
\catcode`\@=11

\magnification 1200
\pretolerance=500 \tolerance=1000 \brokenpenalty=5000

\hsize=134mm \vsize=204mm
\hoffset=0mm \voffset=0mm

\parindent 6mm

\parskip 1,2mm

\baselineskip=12,8pt

\hfuzz 2mm

\newcount\equacount
\def\ifundefined#1{\expandafter\ifx\csname#1\endcsname\relax}
\def\equadef#1{\global\advance\equacount by 1%
  \global\expandafter\edef\csname#1\endcsname{\the\equacount}%
  \the\equacount}
\def\equaref#1{\expandafter\csname#1\endcsname}

\newif\ifpagetitre        \pagetitretrue
\newtoks\hautpagetitre    \hautpagetitre={\hfil}
\newtoks\baspagetitre     \baspagetitre={\hfil}
\newtoks\auteurcourant    \auteurcourant={\hfil}
\newtoks\titrecourant     \titrecourant={\hfil}
\newtoks\hautpagegauche   \newtoks\hautpagedroite
\hautpagegauche={\hfil\tensl\the\auteurcourant\hfil}
\hautpagedroite={\hfil\tensl\the\titrecourant\hfil}

\newtoks\baspagegauche
\baspagegauche={\hfil\tenrm\folio\hfil}
\newtoks\baspagedroite
\baspagedroite={\hfil\tenrm\folio\hfil}

\headline={\ifpagetitre\the\hautpagetitre
\else\ifodd\pageno\the\hautpagedroite
\else\the\hautpagegauche\fi\fi}

\footline={\ifpagetitre\the\baspagetitre
\global\pagetitrefalse
\else\ifodd\pageno\the\baspagedroite
\else\the\baspagegauche\fi\fi}

\font\bbb =msbm10
\def \rr {{\hbox {\bbb R}}}

\font\bbbb =msbm7
\def \rrr {{\hbox {\bbbb R}}}

\font\bbb =msbm10

\font\bbbb =msbm7

\font\bbb =msbm10

\font\bbbb =msbm7

\font\bbb =msbm10

\font\bbbb =msbm7

\font\bbb =msbm10

\font\bbb =msbm10

\font\bbb =msbm10

\font\bbbb =msbm7

\def\ref#1|#2|
       {{\leftskip=6mm\noindent
           \hangindent=5mm\hangafter=1
           \llap{[#1]}\hskip 5mm{#2}.\par}}

\def\pe#1#2 {\pc#1#2|}
\def\pc#1#2|{{\tenrm#1\sevenrm#2}}
\def\pd#1 {{\pc#1|}\ }
\def\bfpc#1#2|{{\tenbf#1\sevenbf#2}}
\def\bfpd#1 {{\bfpc#1|}\ }

\def\cqfd{\hfill\penalty 500\raise 0pt\hbox{\vrule\vbox to 6.5pt{
\hrule width 5.8pt\vfill\hrule}\vrule}\par}

\def \P {{(P_t)}_{t\geq 0}} 
\def \eu {{\rm e}}
\def \L {{\rm L}}
\def \pt {{{(P_t)}_{t \geq 0}}}


\centerline { }

\vskip 10mm

\centerline { \bf ON HARNACK INEQUALITIES}

\vskip 1.5mm

\centerline { \bf  AND OPTIMAL TRANSPORTATION}

\vskip 8mm

\font\pc=cmcsc10 \rm 

\centerline {{\pc D. Bakry, I. Gentil, M. Ledoux}}

\vskip 2mm

\centerline {\it University of Toulouse and University of Lyon, France}

\vskip 1cm

\parindent 1,2cm

{\narrower 
 \parindent 6mm

Abstract. -- {\it We develop connections between Harnack inequalities for the heat flow
of diffusion operators with curvature bounded from below
and optimal transportation. Through heat kernel inequalities,
a new isoperimetric-type Harnack inequality is emphasized.
Commutation properties between the heat and Hopf-Lax
semigroups are developed consequently, providing direct access
to heat flow contraction in Wasserstein spaces.}
\par }

\parindent 6mm

\parskip 1,8mm

\vskip 8mm \goodbreak

{\bf 1. Introduction}

\vskip 3mm

Harnack inequalities classically provide strong tools towards regularity properties 
of solutions of partial differential equations and heat kernel
bounds. A renowned result on the topic is the parabolic inequality by
P. Li and S.-T. Yau [L-Y]
$$ {|\nabla P_t f|^2 \over (P_t f)^2} - {\Delta P_t f \over P_t f} \leq {n \over 2t}  \leqno (1) $$
for the heat semigroup $\pt $ on an $n$-dimensional Riemannian manifold $(M,g)$ with
non-negative Ricci curvature, and every $t>0$ and positive
(measurable) function $f : M \to \rr$. By integration
along geodesics, it yields the Harnack inequality
$$ P_t f(x) \leq P_{t+s} f(y) \Big ( {t+s \over t} \Big )^{n/2} \eu^{d(x,y)^2/4s} \leqno (2) $$
for $f : M \to \rr$ non-negative and $t,s>0$, where $  d(x,y)$ is the Riemannian distance
between $x , y \in M$. The results (1) and (2) admit versions for any lower bound
on the Ricci curvature (cf. [L-Y], [D]). A heat flow proof of (1), in the spirit of the arguments
developed in this work, has been provided in [B-L2].

In the context of diffusion operators, the Harnack inequality (2) may actually loose its relevance due to
the infinite-dimensional feature of some models.
Let $\L = \Delta - \nabla V \cdot \nabla $ be a diffusion operator on
a smooth complete connected Riemannian manifold $(M,g)$,
where $V : M \to \rr$ is a smooth potential, with associated
Markov semigroup $\pt $ and invariant and symmetric measure
$d \mu  = \eu^{-V} dx$ (where $dx$ is the Riemannian volume element).
A notion of curvature-dimension $CD(K, N)$, $K \in \rr$, $N\geq 1$,
of such operators $\L$ has been introduced
by D. Bakry and M. \'Emery [B-\'E] (cf. [B], [Ba-G-L]), which is by now classically
referred to as the $\Gamma _2$ criterion, through the Bochner-type inequality
$$ {1\over 2} \, \L \big (|\nabla f|^2 \big ) - \nabla f \cdot \nabla \L f
    \geq K  |\nabla f|^2 + {1\over N} \, (\L f) ^2 \leqno (3)$$
for any smooth $f : M \to \rr$. (The $\Gamma _2$ operator in this context is precisely
the expression on the right-hand side of (3).) For example, by the standard Bochner formula
from Riemannian geometry, the Laplace operator $\Delta$ on an $n$-dimensional
Riemannian manifold with Ricci curvature bounded from below by $K$ satisfies
the curvature-dimension condition
$CD(K, N)$ with $N \geq n$. On the other hand, on $M = \rr^n$ with
$V$ the quadratic potential, the associated Ornstein-Uhlenbeck operator $\L$ is
intrinsically of infinite dimension $N = \infty$ since (3) cannot hold for some $K \in \rr $ with $N$ finite.
(It actually holds in this example with $K  = 1$ and $N = \infty$.)
In particular a Harnack inequality (2) cannot hold in this case, as well as in further
similar infinite-dimensional models.
Note that, when $N = \infty$, again by the Bochner formula,
the curvature condition $CD(K ,\infty)$, $K \in \rr$, amounts to the local geometric lower bound
$$ {\rm Ric} + {\rm Hess}(V) \geq K  \leqno (4)$$
(as symmetric matrices) uniformly over the manifold (cf. e.g. [Ba-G-L]).

To circumvent the drawbacks attached to the case $K = \infty$,
F.-Y. Wang introduced in [W1] (see also [W2]) a new
form of Harnack inequalities for infinite-dimensional diffusion operators
of the type $\L = \Delta - \nabla V \cdot \nabla $ (and more general ones).
Wang's Harnack inequalities indicate that, under the curvature
condition $CD(K  ,\infty)$ (equivalent to (4)), for every non-negative
(Borel measurable)
function $f $ on $M$, every $t >0$, every $\alpha >1$, and every $x, y \in M$,
$$ \big (P_t f(x) \big )^\alpha  \leq  P_t (f^\alpha ) (y) 
   \, \eu ^{ \alpha d(x,y)^2/2(\alpha -1)\sigma (t)}   \leqno (5) $$
where $\sigma (t) = {1\over K  } \,  (\eu^{2K  t} -1 )$ ($ = 2t$ if $K  = 0$).
The proof of (5) is based on the interpolation
$$  P_s \big ( (P_{t-s} f)^\alpha \big) (x_s), \quad s \in [0,t], $$
along a geodesic ${(x_s)}_{s \in [0,t]}$ joining $x$ to $y$
together with the commutation, for all $t \geq 0$ and smooth $g: M \to \rr$,
$$ | \nabla P_t g| \leq \eu^{-K  t} P_t \big ( |\nabla g| \big ) \leqno (6)$$
as an equivalent formulation of the curvature lower bound $CD(K  ,\infty)$
(cf. [Ba-G-L]). In a sense, the gradient bound (6) may be thought of as the counterpart
of the Li-Yau inequality (1) in this context.

Changing $f$ into $f^{1/\alpha }$, in the asymptotics
$f ^{1/\alpha } \sim 1 + {1 \over \alpha }  \log f $ as $\alpha \to \infty$
in (5), a log-Harnack inequality
$$ P_t (\log f) (x) \leq \log P_t f (y) + {d(x,y)^2 \over 2\sigma (t)}  \leqno (7)$$
also holds (cf. [Bo-G-L], [W3]).
It was further shown in [W1], [W3] that either (5) (for one
$\alpha >1$) or (7) imply back, as $ t \to 0$, the curvature condition
$CD(K  ,\infty)$ (in its infinitesimal form (4)).

The aim of this work is two-fold. We will first show how the previous 
infinite-dimensional Harnack inequalities
may actually be seen as consequences of a suitable
functional inequality of isoperimetric type. On the basis of this observation, we establish next
a kind of isoperimetric-type Harnack inequality.
These results naturally lead to develop connections between isoperimetric-type 
Harnack inequalities (in direct or reverse form) and commutation properties between
diffusion and Hopf-Lax semigroups.
By the dual Kantorovich optimal transportation formalism, 
Wasserstein contraction properties along the heat flow are then derived.

Two observations are actually at the starting point of this work.
For simplicity in the (somewhat informal) discussion below, we restrict ourselves to 
the curvature condition $CD (0, \infty)$ (with thus $K  =0$).

First, the gradient bound (6) (and thus the curvature condition $CD (0, \infty)$)
is known to imply (to be equivalent) to
logarithmic Sobolev inequalities under the heat kernel measures
$P_t$, in particular in reverse form
$$ t \, {|\nabla P_t f|^2 \over P_t f} \leq  P_t (f \log f ) - P_t f \log P_t f \leqno (8) $$
for every (bounded measurable) $f >0$ and every $t>0$ (cf. [Ba-G-L]). Inequalities like the preceding one
are understood point-wise throughout this work.
Now, as was noticed by M.~Hino~[H], the latter ensures that whenever $0 < f \leq 1$ and
$ \psi  = \sqrt { \log (1/P_t f)}$, then
$$  | \nabla \psi  |^2 \leq {1 \over 2t} \, .$$
In other words, $\psi $ is Lipschitz with Lipschitz coefficient
less than or equal to $(2t)^{-1/2}$. In particular, for every $x, y \in M$,
$$  \sqrt  {\log {1 \over P_t f(x)} }
      \leq  \sqrt {\log {1 \over P_t f(y)} }  + { d (x,y)\over {\sqrt {2t}}}  $$
where we recall that $  d(x,y)$ is the Riemannian distance between $x$ and $y$.
After some work, it may then
be shown that for each $\varepsilon >0$, there exists $C(\varepsilon ) >0$ such that
$$ \big ( P_t f(x) \big )^2 \leq  C(\varepsilon ) \, P_t (f^2) (y) \,
     \eu^{d(x,y)^2/2(1+\varepsilon )t} , $$
that is as close as possible to (5) (for $\alpha =2$). 

It should be mentioned that it is precisely the dimensional version of the reverse logarithmic
Sobolev inequality (8) which has been used in [B-L2] to provide a monotonicity proof of the
Li-Yau parabolic inequality (1). We will exploit this information towards dimensional
statements in Section 5 below. For further dimensional Harnack inequalities under
the curvature-dimension condition $CD(K  ,N)$, comparing in particular different
times, see [W4], [E-K-S], [K2].

The second observation at the starting point of this investigation is the links
between Harnack-type inequalities and optimal transportation already put forward
in [Bo-G-L] where semigroup tools were developed towards a proof
of the Otto-Villani HWI inequality [O-V] (cf. [V1], [V2]). We briefly recall the basic step.
Namely, the log-Harnack inequality (7) may be translated equivalently as
$$ P_t (\log f)  \leq  Q_{2t} (\log P_t f) \leqno (9) $$
where ${(Q_s)}_{s >0}$ is the Hopf-Lax infimum-convolution semigroup
$$ Q_s \varphi  (x) = \inf_{y \in M} \Big [ \varphi (y) + {d(x,y)^2 \over 2s} \Big],
   \quad x \in M, \, \, s >0. $$
Assume now that $\mu $ is a probability measure 
and let $f >0$ be a (bounded) probability density with respect to $\mu $. Then,
by time reversibility and (9) applied to $P_t f$, $ t >0$,
$$ \int_M P_t f \log P_t f d\mu  = \int_M f P_t (\log P_t f) d\mu 
    \leq \int_M f \, Q_{2t} (\log P_{2t} f) d\mu  . $$
Now $\int_M \log P_{2t} f d\mu  \leq 0 $ by Jensen's
inequality. Hence, combining with the scaling properties of ${(Q_s)}_{s >0}$,
$$ \int_M P_t f \log P_t f d\mu  
      \leq  {1 \over 2t} \bigg [ \int_M Q_1 \varphi  \, f d\mu  - \int_M \varphi  d\mu  \bigg ]$$
where $ \varphi  = 2t \log P_{2t} f$.
Recall then the (quadratic) Wasserstein distance $ {\rm W}_2(\nu ,\mu )$ 
between two probability measures $\mu $ and $\nu$ on $M$ defined by
$$ {\rm W}_2(\nu ,\mu ) = \bigg ( \int_{M \times M} d(x,y)^2 d\pi (x,y) \bigg)^{1/2}$$
where the infimum is taken over all couplings $\pi $ with respective marginals $\nu $ and $\mu $.
The Kantorovich dual description
$$ {1 \over 2} \, {\rm W}_2(\nu ,\mu )^2 
    = \sup \bigg ( \int_M  Q_1 \varphi  \,  d\nu - \int _M \varphi   d\mu  \bigg )  \leqno (10) $$
where the supremum runs over all bounded continuous functions $\varphi  : M \to \rr$
(cf. e.g. [V1]) then yields with $d\nu = f d\mu $ 
$$ \int_M P_t f \log P_t f d\mu  \leq  {1 \over 4t} \, {\rm W}_2^2 (\nu , \mu ). \leqno (11) $$

Note that the preceding argument similarly yields, for every $ t > 0$,
$$  \int_M P_t f \log P_t f d\mu  \leq  {1 \over 4 t} \, {\rm W}_2^2 (f\mu , g\mu )
    +  \int_M g \log g \, d\mu \leqno (12) $$
where $g$ is a further probability density with respect to $\mu $, and where,
for simplicity here, $f\mu $ and $g\mu $ denote the probability measures $fd\mu $ and
$gd\mu $. Indeed, write in the preceding notation that
$$  \int_M P_t f \log P_t f d\mu  
      \leq  {1 \over 2t} \bigg [ \int_M Q_1 \varphi  \, f d\mu  - \int_M \varphi  g d\mu  \bigg ]
         +  \int_M g \log P_{2t} f d\mu .$$
Since by convexity $ \int_M g \log P_{2t} f d\mu \leq \int_M g \log g d\mu $, the claim follows.

Now (11) is actually the major step in the semigroup proof of the HWI inequality
of [O-V] under the curvature condition $CD (0, \infty)$.
Namely, the classical heat flow interpolation scheme
(cf. [B], [Ba-G-L]) indicates that, for every suitable 
smooth probability density $f : M \to \rr$ and every $t\geq 0$,
$$ \eqalign {
    \int_M  f \log  f d\mu   -  \int_M P_t f \log P_t f d\mu  
    & = - \int_0^t \bigg ({d \over ds} \int_M P_s f \log P_s f d\mu \bigg) ds \cr
     &  = \int_0^t   \int_M {|\nabla P_s f|^2 \over P_s f } \, d\mu  \, ds. \cr } $$
Since $|\nabla P_s f| \leq P_s (|\nabla f|)$ according to (6),
$$   {|\nabla P_s f|^2 \over P_s f } \leq  P_s \Big ( {|\nabla f|^2 \over f } \Big )  $$
by the Cauchy-Schwarz inequality along the Markov kernel $P_s$. Therefore,
$$ \int_M  f \log  f d\mu   \leq  \int_M P_t f \log P_t f d\mu   
      + t  \int_M {|\nabla f|^2 \over f } \, d\mu  .$$
Together thus with (11), optimization in $t >0$ yields
$$ \int_M  f \log  f d\mu 
    \leq  {\rm W}_2 (\nu , \mu ) \bigg (  \int_M {|\nabla f|^2 \over f } \, d\mu  \bigg)^{1/2}  $$
which is the announced HWI inequality, connecting Entropy, Wasserstein distance and Fisher Information.
Similar arguments may be developed under $CD (K  , \infty)$ for any
$K  \in \rr$ to yield the full formulation of Otto-Villani's HWI inequality
(cf. [Bo-G-L], [Ba-G-L]). Note that together with (12), the argument also recovers the
known inequality (cf. e.g. [CE])
$$ \int_M  f \log  f d\mu 
    \leq  {\rm W}_2 (f\mu , g\mu ) \bigg ( \int_M {|\nabla f|^2 \over f } \, d\mu  \bigg)^{1/2}  
       + \int_M  g \log  g d\mu $$
for probability densities $f$ and $g$ with respect to $\mu$.

For the matter of comparison, it might be worthwhile mentioning that the recent
Kuwada lemma (see [G-K-O]) develops similar arguments towards the inequality
$$ W_2^2 (P_tf \mu , f\mu ) 
     \leq t  \bigg [ \int_E f \log f d\mu - \int_E P_t f \log P_t f d\mu \bigg]  \leqno (13)$$
for any probability density $f$ with respect to $\mu $ and any $t \geq 0$.
Indeed, for $ \varphi : E \to \rr$ bounded and continuous,
$$ \eqalign {
   \int_M Q_1 \varphi \, P_t f d\mu  - \int_M \varphi f d\mu 
      & = \int_0^1 \bigg ({d \over ds} \int_M  Q_s \varphi P_{st} f d\mu  \bigg) ds \cr
      & = \int_0^1 \int_M \Big [ - {1\over 2} \, |\nabla Q_s \varphi |^2 P_{st} f 
                + t Q_s \varphi \, \L P_{ts} f \Big ] d\mu \, ds \cr
         & = \int_0^1 \int_M \Big [ - {1\over 2} \, |\nabla Q_s \varphi |^2 P_{st} f 
                -  t \, \nabla Q_s \varphi \cdot \nabla P_{st} f \Big ] d\mu \, ds  \cr} $$
by the fact that the Hopf-Lax semigroup solves the standard Hamilton-Jacobi equation
and by integration by parts.
Next, by the Cauchy-Schwarz inequality,
$$     \int_M Q_1 \varphi \, P_t f d\mu  - \int_M \varphi f d\mu 
    \leq  {t^2 \over 2} \int_0^1 \int _M {|\nabla P_{st} f|^2 \over P_{st} f} \, d\mu  \, ds $$
which yields (13) by the Kantorovich duality and integration.

The inequality (13) is actually at the core of the gradient flow interpretation of the
heat flow in Wasserstein space (cf. [J-K-O], [O], [V1], [V2], [A-G-S1]), and immediately follows
for example from
the Benamou-Brenier [B-B] dynamical characterization of the Wasserstein distance
in smooth spaces such as Riemannian manifolds. Kuwada's argument above extends its
validity to a general, possibly nonlinear, setting.

From the point of view of functional inequalities, 
(13) somewhat works in the other direction with respect to (11). 
Namely, while (11) leads to the HWI inequality, (13)  has been identified
in [G-L] at the root of the Otto-Villani theorem [O-V] (cf. [Bo-G-L], [V1], [V2])
connecting logarithmic Sobolev inequalities to transportation cost inequalities.

\medskip

On the basis of these two early observations, the
purpose of this work is, as announced, to develop a synthetic and refined treatment of
Harnack-type inequalities for diffusion operators with curvature bounded from below
and of their connections
with transportation cost inequalities. The various contributions of this work
are summarized as follows.

In Section 2, we provide a direct treatment of Wang's Harnack inequalities
(5) and (7) relying on an improved,
isoperimetric-type version, of the reverse logarithmic Sobolev inequality (8) along the heat flow. 

This reverse isoperimetric-type inequality in turn
implies a new isoperimetric version of Harnack inequalities emphasized in Section 3.
For example, under non-negative curvature $CD(0, \infty)$, it yields that
for any (Borel) measurable set $A$ in $M$, any $t>0$ and any $x, y \in M$,
$$ P_t({\bf 1}_A) (x) \leq  P_t \big ({\bf 1}_{A_{d(x,y)}} \big ) (y) \leqno (14) $$
where $A_\varepsilon $ is the $\varepsilon $-neighborhood of
$A$ in the metric $d$. This result seems to be new even for the standard heat flow operator
on a Riemannian manifold. It is optimal for the standard heat kernel on $\rr^n$
as is immediately checked on the explicit kernel representation.

A direct consequence of (14) is the commutation
$$ P_t (Q_s) \leq Q_s(P_t), \quad t,s >0, \leqno (15) $$
between the heat and Hopf-Lax semigroups (under non-negative curvature)
which we emphasize in Section 4. This commutation was actually used
earlier by K. Kuwada [K1] in the study of the duality of gradient estimates
at the root of the contraction property of the Wasserstein distance along the heat flow 
$$  {\rm W}_2 (\mu_t , \nu_t )  \leq  {\rm W}_2 (\mu_0 , \nu_0 )  $$
where $d\mu_t = P_t f d\mu $ and $ d\nu _t = P_t g d\mu $, $t \geq 0$,
$f,g$ probability densities with respect to $\mu $.
Such contraction properties have been investigated in this context 
in [O], [C-MC-V], [vR-S], [O-W] (see also [W2], [V1], [V2]), and are also,
following [A-G-S1], [Er], a main consequence of the EVI approach
discussed in Section 6.

The commutation property (15)
may actually be reached in several ways, and Section~5
presents a variety of methods depending on the underlying context, including
the original approach of [K1]. This section further includes dimensional versions
of the commutation property together with the corresponding Wassertein contractions.

In the last Section 6, we briefly discuss some connections between the material
presented here and recent developments, following [A-G-S3], around
the Evolutionary Variational Inequality (EVI) expressing in the preceding notation that
$$ {\rm W}_2^2 (\mu _t , \nu _0 ) +  2t \int_M P_t f \log P_t f d\mu 
    \leq  {\rm W}_2^2 (\mu _0 , \nu _0 ) + 2 t \int_M g \log g d \mu .$$
This property actually connects the $\Gamma _2$ Bakry-\'Emery $CD (K  , \infty)$ curvature
condition ([B-\'E], [B], [Ba-G-L]), expressed by the commutation (6),
with the curvature bound in the
sense of Lott-Villani-Sturm in metric measure spaces as convexity of relative
entropy along the geodesics of optimal transportation ([L-V], [S1], [S2], [V2]).
The recent main achievement by L. Ambrosio, N. Gigli and G. Savar\'e [A-G-S2], [A-G-S3] actually provides 
a link between the $\Gamma _2$ and Lott-Villani-Sturm curvature lower bounds in the class
of the Riemannian energy measure spaces through the EVI. In Section 6, we sketch,
following [A-G-S3],
the principle of proof of the EVI in a smooth setting, for comparison with some
of the tools developed here.

\vskip 2mm

For simplicity in the exposition, the results of this work are presented
in the weighted Riemannian setting, for thus diffusion operators $\L = \Delta - \nabla \cdot \nabla V$
on a complete connected Riemannian manifold $(M,g)$ with invariant 
and reversible measure $d\mu = \eu^{-V} dx$ (not necessarily a probability measure)
where $V : M \to \rr$ is a smooth potential. Integration by parts with respect to $\L$
is expressed by $\int_M f (-\L g) d\mu  = \int_M \nabla f \cdot \nabla g d\mu $ for
smooth functions $f, g : M \to \rr$. The associated curvature condition
$CD (K  , \infty)$, $K  \in \rr$, is expressed equivalently by 
(4) as the infinitesimal version of the Bochner-type inequality (3). It amounts
to the standard lower bound on the Ricci curvature for the Laplace operator $\Delta $ on $(M,g)$. 
The curvature condition $CD(K  ,\infty)$ is also equivalent to
the gradient bound (6) which is, in an essential manner, the only way the curvature
condition will be used throughout this work.

We refer to the general references [B], [Ba-G-L] for a precise description
of this framework and the relevant properties. Most of the results below actually extend
to the more general setting of a Markov diffusion Triple $(E, \mu, \Gamma )$ emphasized in
[B], [Ba-G-L], consisting of a state space $E$
equipped with a diffusion semigroup $\P$ with infinitesimal generator $\L$
and carr\'e du champ operator $\Gamma $ and
invariant and reversible $\sigma $-finite measure $\mu $. In the weighted
Riemannian context, $ \Gamma (f,f) = |\nabla f|^2$ for smooth functions.
In this setting, the abstract curvature condition $CD(K  , \infty)$, $K  \in \rr$, stems
from the Bochner-type inequality (3) (with $N = \infty$) and the abstract
$\Gamma _2$ operator going back to [B-\'E] (see [B], [Ba-G-L]).
The condition $CD(K  , \infty)$ is equivalent to the gradient bound (6)
$$ \sqrt {\Gamma (P_t f)} \leq \eu^{-K  t} P_t  \big (\sqrt {\Gamma (f)} \, \big )  $$
for every $t\geq 0$ and every $f$ in a suitable algebra of functions.
The state space $E$ may be endowed with an intrinsic distance $d$ for which Lipschitz
functions $f$ are such that $\Gamma (f) $ is bounded ($\mu $-almost everywhere).
Note also that at the level of the local inequalities along the semigroup, generators of the type
$\L = \Delta + Z$ for some smooth vector field $Z$ on a manifold $M$
may be covered similarly as developed in [W1], [W2], [W3].

\vskip 8mm \goodbreak

{\bf 2. Reverse isoperimetry and Wang's Harnack inequalities}

\vskip 3mm

In this section, we address a direct proof of Wang's Harnack inequalities (5) and (7) along
the Hino argument on the basis of a reinforced family of heat kernel inequalities
first emphasized in [B-L1]. 

Denote by $I : [0,1] \to \rr_+$ the Gaussian isoperimetric function defined by
$ I = \varphi \circ \Phi ^{-1}$ where
$$ \Phi (x) = \int_{-\infty}^x \eu^{-u^2/2} {du \over \sqrt {2\pi }}, \quad x \in \rr,$$
and $ \varphi = \Phi '$. The function $I$ is concave continuous, symmetric with respect to the vertical
line going through ${1 \over 2}$ and such that $I(0) = I(1) = 0$,
and satisfies the basic differential equality $I \, I'' = -1$.
Moreover $ I (v) \sim v \sqrt {2 \log {1\over v} } $ as $ v \to 0$.

The following statement, as a kind of reverse isoperimetric-type
inequality, was first put forward in [B-L1]. We enclose a proof for completeness
(see also [Ba-G-L]).

\vskip 2mm

{\bf Proposition 2.1.} {\sl Under the curvature condition $CD(K  ,\infty)$ for some $K  \in \rr$,
for every (measurable) function $f$ on $M$ with values in $[0,1]$ and every $t > 0$,
$$ \big [I (P_tf)\big ]^2- \big [P_t\big (I  (f) \big )\big ]^2
		\geq   \sigma (t)  |\nabla P_t f|^2    \leqno (16)$$
where $\sigma (t) = {1\over K  } \,  (\eu^{2K  t} -1 )$ ($ = 2t$ if $K  = 0$).}

\vskip 2mm

{\it Proof.} Work with a function $f$ such that $  \varepsilon  \leq f \leq 1 - \varepsilon $
for some $\varepsilon > 0$.
By the heat flow interpolation, write
$$ \big [I (P_tf)\big ]^2- \big [P_t\big (I  (f) \big )\big ]^2
	= - \int _0^t {d\over ds }\big [P_s\big (I (P_{t-s}f)\big )\big ]^2 ds.$$
Now, by the chain rule for the diffusion operator $\L$,
$$ \eqalign { - {d\over ds }\big [P_s\big (I (P_{t-s}f) \big )\big ]^2
& = - 2 P_s\big (I (P_{t-s}f) \big )  P_s \Big ( \L I (P_{t-s} f) - I' ( P_{t-s}f) \, \L P_{t-s}f) \Big) \cr
&=  - 2 P_s\big (I (P_{t-s}f) \big ) 
		P_s \big (I ''(P_{t-s}f) |\nabla P_{t-s} f|^2 \big )\cr
&= 2 P_s\big (I(P_{t-s}f) \bigr )
		 P_s\bigg ( {|\nabla P_{t-s} f|^2 \over I (P_{t-s}f)} \bigg )\cr }$$
where we used that $ I \, I'' = -1$ in the last step.
Since $P_s$ is given by a kernel, it satisfies a Cauchy-Schwarz inequality,
and hence
$$  P_s(Y) P_s\bigg ( {{X^2} \over  Y} \bigg ) \geq  \big [P_s(X)\big ]^2, \quad X,Y \geq 0.$$
Hence, with $X = |\nabla P_{t-s} f|$ and $Y = I (P_{t-s}f)$,
$$ \big [I (P_t f)\big ]^2- \big [P_t\big (I  (f) \big )\big ]^2
		\geq  2 \int _0^t \Big [P_s\big ( |\nabla P_{t-s} f| \,  \big )\Big ]^2 ds. $$
By the gradient bound (6) applied to $g = P_{t-s}f$, it follows that
$$ \big [I (P_tf)\big ]^2- \big [P_t\big (I  (f) \big )\big ]^2
		\geq  2 \int _0^t \eu^{2K  s} ds \, |\nabla P_t f|^2 $$
which is the result. \cqfd			

\vskip 2mm

For the comparison with the Hino observation mentioned in the introduction, note that
(16) of Proposition 2.1 implies the reverse logarithmic Sobolev inequality (8)
by applying it to $\varepsilon f$ and letting $\varepsilon \to 0$.

As announced, we next show how Proposition 2.1 actually covers Wang's Harnack inequalities
recalled in the introduction. The main consequence is put forward
in the following corollary that actually entails most of the inequalities emphasized in this work.

\vskip 2mm

{\bf Corollary 2.2.} {\sl Under the curvature condition $CD(K  ,\infty)$ for some $K  \in \rr$,
for every (measurable) function $f$ on $M$ with values in $[0,1]$, every $t > 0$
and every $x, y \in M$,  
$$ \Phi^{-1} \circ P_t f (x) \leq \Phi^{-1} \circ P_t f (y) + {d (x,y)\over \sqrt {\sigma (t)}}  \leqno (17) $$
where $ d(x,y)$ the Riemannian distance between $x$ and $y$.}

\vskip 2mm

Namely, in terms of gradient bounds, Proposition 2.1 implies that for every $f$ with
values in $[0,1]$,
$$ |\nabla P_t f| \leq  {1 \over \sqrt {\sigma (t)}} \, I (P_tf).  $$
Since $ (\Phi^{-1})' = {1 \over I}$, it follows that $\Phi^{-1} \circ P_t f $ is
$\sigma (t)^{-1/2}$-Lipschitz, $t>0$, which amounts to the corollary.

Towards a first illustration of Corollary 2.2,
set $ \delta = d(x,y)/\sqrt {\sigma (t)}$, so that (17) reads as
$$ P_t f(x) \leq \Phi \big ( \Phi^{-1} \circ P_t f (y) + \delta \big) .$$
Apply now this inequality to
${\bf 1}_{\{ f \geq a\}}$, $a \geq 0$, for a non-negative function $f$ on $M$.
Denoting by $\lambda  $ the distribution of $f$ under $P_t$ at the point $y$
(that is $\lambda (B) = P_t ({\bf 1}_{\{ f \in B\}} ) (y)$ for every Borel set $B$ in $\rr$),
$$ P_t ({\bf 1}_{\{ f \geq a\}})(x) 
     \leq \Phi \Big ( \Phi^{-1} \big(\lambda  ([a, \infty)) \big) + \delta \Big) .$$
Integrating in $a \geq 0$ and using Fubini's theorem, denoting by
$d\gamma (u) = \eu^{-u^2/2} {du \over \sqrt {2\pi }}$ the standard
Gaussian distribution on the line,
$$ P_tf(x) \leq \int_0^\infty \! \! \int_{-\infty}^{\Phi^{-1} (\lambda  ([a, \infty))) + \delta} d\gamma (u) da
       = \int_{- \infty}^\infty  \bigg (
         \int_0^\infty {\bf 1}_{\{ u \leq \Phi^{-1} (\lambda  ([a, \infty))) + \delta \} } da \bigg) d\gamma (u) . $$
Change $u$ into $u + \delta$ to get
$$ P_tf(x)  \leq \eu^{-\delta ^2/2} \int_{-\infty}^\infty \eu^{-\delta u}
            \bigg ( \int_0^\infty {\bf 1}_{\{ \Phi (u) \leq  \lambda  ([a, \infty)) \} } da \bigg) d\gamma (u) .$$
Changing $u$ into $-u$ and denoting by $F$ the distribution function of $\lambda  $, it follows that
$$ P_tf(x)  \leq \eu^{-\delta ^2/2} \int_{-\infty}^\infty \eu^{\delta u}
            \bigg ( \int_0^\infty {\bf 1}_{\{ F(a) \leq \Phi (u) \} } da \bigg) d\gamma (u) .$$            
After the further change of variables $ v = \Phi (u)$,
$$ P_tf(x)  \leq \eu^{-\delta ^2/2} \int_0^1 \eu^{\delta \, \Phi^{-1}(v)}
            \bigg ( \int_0^\infty {\bf 1}_{\{ F(a) \leq v \} } da \bigg) dv .$$

The next statement summarizes the conclusion reached so far.

\vskip 2mm

{\bf Theorem 2.3.} {\sl Under the curvature condition $CD(K  ,\infty)$ for some $K  \in \rr$,
for every non-negative (measurable) function $f$ on $M$, every $t > 0$ and every $x, y \in M$,
$$ P_tf(x)  \leq \eu^{-\delta ^2/2} \int_0^\infty \eu^{\delta  \,\Phi^{-1}\circ F(r)} \, r\, dF(r) $$
where $\delta = d(x,y)/\sqrt {\sigma (t)}$ and
$F (r) = P_t ( {\bf 1} _{\{ f \leq r \} } ) (y)$, $ r \geq 0$,
is the distribution function of $f$ under $P_t$ at the point $y$.}

\vskip 2mm

Theorem 2.3 appears at the root of the various Harnack inequalities in this context.
It is however not expressed in a very tractable form. But it easily implies known ones. For example,
by Cauchy-Schwarz,
$$ \eqalign {
   \int_0^\infty \eu^{\delta \, \Phi^{-1}\circ F(r)} \, r \, dF(r)
   & \leq \bigg ( \int_0^\infty \eu^{2\delta \, \Phi^{-1}\circ F(r)} dF(r) \bigg)^{1/2}
          \bigg ( \int_0^\infty  r^2 \, dF(r) \bigg )^{1/2} \cr
     & \leq  \eu^{\delta^2} \big (P_t(f^2)(y) \big)^{1/2} \cr } $$
since
$$ \int_0^\infty \eu^{2\delta \Phi^{-1}\circ F(r)} dF(r)
    = \int_0^1 \eu^{2 \delta \Phi^{-1} (v)} dv = \int_{-\infty}^\infty \eu^{2\delta u} d\gamma (u)
        = \eu^{2\delta^2}. $$     
The preceding therefore yields Wang's Harnack inequality (5) for $\alpha =2$,
$$ P_tf (x)^2  \leq  P_t(f^2)(y) \, \eu^{d(x,y)^2/\sigma (t)} . $$
By H\"older's inequality rather than Cauchy-Schwarz, one obtains the whole family of inequalities (5) with
$\alpha >1$. Using the entropic inequality yields similarly the log-Harnack inequality (7).
With respect to Wang's original argument,
the proof here avoids interpolation along geodesics (although the length space property
is required to move from (16) to (17)).

\vskip 8mm \goodbreak

{\bf 3. Isoperimetric-type Harnack inequalities}

\vskip 3mm

As announced, the basic Lipschitz inequality (17) of Corollary 2.2
may be seen at the origin of a number of
inequalities of interest, and this section develops further consequences 
in combination with isoperimetric bounds for heat kernel measures.
To this task, recall first the isoperimetric comparison theorem for heat kernel measures
under curvature bounds of [B-L1]. Recall $I$ the Gaussian isoperimetric function.

\vskip 2mm

{\bf Theorem 3.1.} {\sl Under the curvature condition $CD(K  ,\infty)$ for some $K  \in \rr$,
for every smooth function $f$ on $M$ with values in $[0,1]$ and every $t\geq 0$,
$$ I  (P_tf)  \leq  P_t \Big (\sqrt {I ^2(f) + K   (t)|\nabla f|^2 }\, \Big ) $$
where $K   (t) = {1\over K  } \,  (1 - \eu^{-2K  t}  )$ ($ = 2t$ if $K  = 0$).}

\vskip 2mm

As developed in [B-L1] (cf. also [Ba-G-L]), this result is an isoperimetric comparison theorem
expressing that the isoperimetric profile of the heat kernel measures is bounded from below
by the Gaussian isoperimetric function (up to a scaling depending on $t$ and $K  $).
This comparison may classically (cf. [B-H], [B-L1]) 
be translated in terms of isoperimetric neighborhoods in the sense that
for any Borel measurable (or closed) set $A \subset M$ and any $\varepsilon >0$,
$$ P_t ({\bf 1}_{A_\varepsilon }) (y)  
     \geq   \Phi \bigg ( \Phi^{-1} \big (P_t ({\bf 1}_A) (y) \big) 
           + {\varepsilon \over \sqrt {K  (t)}} \bigg)  \leqno (18)$$
where $A_\varepsilon $ is the (open) $\varepsilon $-neighborhood of $A$ in the distance $d$,
for any $y \in M$ and $t>0$.

Applied to $f = {\bf 1}_A$, the Lipschitz property (17) ensures on the other hand that,
for any measurable set $ A \subset M$, and again with $\delta = d (x,y) / \sqrt {\sigma (t)}$,
$$ P_t ({\bf 1}_A) (x) \leq   \Phi \Big ( \Phi^{-1} \big (P_t ({\bf 1}_A) (y) \big) + \delta \Big ) .
    \leqno (19) $$   
The combination of (18) and (19) together with the fact that
${K  (t) \over  \sigma (t) } = \eu^{-2K  t}$ then yields the following isoperimetric-type
Harnack inequality.

\vskip 2mm

{\bf Theorem 3.2.} {\sl Under the curvature condition $CD(K  ,\infty)$ for some $K  \in \rr$,
for every measurable set $A$ in $M$, every $t\geq 0$ and every $x, y \in M$
such that $d (x,y) >0$,
$$ P_t ({\bf 1}_A) (x)  \leq  P_t \big ({\bf 1}_{A_{d_t}} \big ) (y)  $$
where $ d_ t = \eu^{-K  t}d (x,y)  $. In particular, when $K  =0$,}
$$ P_t ({\bf 1}_A) (x)  \leq  P_t \big ({\bf 1}_{A_ {d(x,y)}}\big ) (y) .  $$

\vskip 2mm

\vskip 8mm \goodbreak

{\bf 4. The commutation property and contraction in Wasserstein distance}

\vskip 3mm

The isoperimetric-type Harnack inequality of Theorem 3.2 has several consequences of interest in terms
of commutation properties between the heat and the Hopf-Lax semigroups  which
in turn entails the contraction property of the heat flow with respect to Wassertein metrics.

Recall the Hopf-Lax infimum-convolution semigroup (cf. [E], [V1], [V2])
$$ Q_s f(x) = \inf_{y \in M} \Big [ f(y) + {d(x,y)^2 \over 2s} \Big] , \quad
    x \in M, \, \, s >0. $$
The announced commutation
property was actually used first by K. Kuwada [K1] in the analysis of gradient bounds and
Wasserstein contractions (see Corollary 4.2 below). The proof in [K1], developed in the context of length
spaces and actually for more general costs, relies on an interpolation along geodesics and the use of the
Hamilton-Jacobi equation (see the first alternate proof in Section 5 below).

\vskip 2mm

{\bf Theorem 4.1.} {\sl Under the curvature condition $CD(K  ,\infty)$ for some $K  \in \rr$,
for any $t, s >0$ and any bounded continuous function $f: M \to \rr$,}
$$ P_t(Q_s f) \leq Q_{\eu^{2K  t}s} (P_t f) . \leqno (20)$$

\vskip 2mm 

{\it Proof.} Let without loss of generality $f$ be non-negative on $M$.
It is enough by homogeneity to consider $s=1$.
Let $x, y $ be arbitrary (distinct) fixed points in $ M$ and set $d_t = \eu^{-K  t} d(x,y) > 0$.
Set $ A = \{ Q_1 f \geq a\}$ for $ a \geq 0$.
If $ z \in  A_{d_t} $, there exists $\xi \in A$ such that $ d(z, \xi ) \leq d_t$ so that
$$ f(z) + {d_t^2 \over 2} \geq f(z) + {d(z, \xi )^2 \over 2}  \geq Q_1 f(\xi ) \geq a . $$
Hence $A_{d_t} \subset \{ f + d_t^2/2 \geq a\}$. Therefore, by Theorem 3.2,
$$ P_t ({\bf 1}_{\{Q_1f \geq a\}}) (x)  \leq  P_t ({\bf 1}_{\{f + d_t^2/2 \geq a\} }) (y) . $$
Integrating in $a \geq 0$ yields
$$ P_t(Q_1 f)(x)  \leq P_tf(y) + {d_t^2\over 2} \, .$$
Taking then the infimum in $y \in M$ yields the result by
definition of the infimum-convolution $Q_1$. \cqfd

\vskip 2mm

The infimum-convolution semigroup
${(Q_s)}_{s>0}$ being solution of the Hamilton-Jacobi equation
$ \partial_s u = - {1\over 2} \, |\nabla u|^2$ with initial condition $u(0,\cdot ) = f$,
the commutation property (20) implies by a Taylor expansion at $s=0$ that
$ |\nabla P_t f|^2 \leq \eu^{-2K  t} P_t \big (|\nabla f|^2 \big)$ for every $t \geq 0$.
This gradient bound, weaker than (6), is still equivalent to the curvature bound $CD(K  ,\infty)$
(cf. [Ba-G-L]), providing therefore a converse to Theorem 4.1. In particular also, the isoperimetric
Harnack inequality from Theorem 3.2 is actually equivalent to the curvature condition $CD(K  ,\infty)$.

As announced, it immediately follows from the commutation property (20)
of Theorem 4.1 that the Wasserstein
distance $W_2$ is contractive along the semigroup $\pt$, an observation
due to K. Kuwada [K1]. The Wasserstein contraction property in this context
may be traced back in the investigation [O] of the heat flow
as a gradient flow in the Wasserstein space, further developed in [C-MC-V].
Further contributions include [vR-S] with a stochastic proof, [O-W]
with an Eulerian point of view, or [A-G-S1], [Er] in connection
with the EVI (Section 6). See also [V1], [V2].
The proof presented here on the basis of Theorem 4.1 extends to the abstract
Markov semigroup setting of [B], [Ba-G-L]. 
The measure $\mu $ is assumed here to be a probability measure.

\vskip 2mm

{\bf Corollary 4.2.} {\sl Under the curvature condition $CD(K  ,\infty)$ for some $K  \in \rr$,
for any $t \geq 0$,
$$ W_2 (\mu _t , \nu_t ) \leq \eu^{-2K  t} \, W_2 (\mu_0  , \nu_0 )  \leqno (21) $$
where $d \mu _t = P_t f d\mu $ and $d\nu _t = P_t g d\mu $ for
probability densities $f, g $ with respect to
the probability measure $\mu $.}

\vskip 2mm 

{\it Proof.} For any bounded continuous $\varphi : M \to \rr$,
by time reversibility and the commutation property (20),
$$ \eqalign {
    \int_M Q_1 \varphi  \, P_t f d\mu  - \int_M \varphi \, P_t g d\mu 
    & = \int_M P_t (Q_1 \varphi ) f d\mu  - \int_M P_t\varphi \,  g d\mu  \cr
    & \leq  \int_M Q_{\eu^{2K  t}} (P_t \varphi ) f d\mu  - \int_M P_t\varphi \,  g d\mu  \cr
    & \leq \eu^{-2K  t} \bigg [ \int_M Q_1 ( \eu^{2K  t} P_t \varphi ) f d\mu 
            -   \int_M \eu^{2K  t} P_t\varphi \,  g d\mu \bigg ] \cr
    & \leq {\eu^{-2K  t} \over 2} \,   \, W_2^2 (\mu _0 , \nu _0 ) \cr }$$
where the last step follows from the Kantorovich dual description (10) of
the Wasserstein distance $W_2$.  The proof is complete. \cqfd

\vskip 2mm

By adapting Theorem 4.1 to costs $d(x,y)^p$, the
same argument works for any Wasserstein distance $W_p$, $1 \leq p < \infty$,
extending the contraction property of Corollary 4.2 to this class.
More general Wasserstein functionals associated to further transportation costs
may be considered similarly. In [K1], K. Kuwada established
the equivalence of the Wasserstein contraction property for the cost $d(x,y)^p$
with the bound (6) with power $q$ on the gradient, where $q>1$ and $p<\infty$ are
dual exponents.

Note that one main conclusion of the work
[vR-S] by M.-K. von Renesse and K.-T. Sturm is the equivalence of the Wasserstein
contraction of Corollary 4.2 with the curvature bound. Actually, reading backwards the proof of Corollary 4.2,
the contraction property (21) indicates that for all $\varphi : M \to \rr$
bounded and continuous,
$$ \int_M P_t (Q_1 \varphi )  f d\mu  - \int_M P_t\varphi \, g d\mu 
   \leq  \int_M Q_1 \varphi  P_t f d\mu  - \int_M \varphi P_t g d\mu 
     \leq { \eu^{-2K  t} \over 2}  \, W_2^2 (\mu _0 , \nu _0 )  . $$
Now, if $x, y \in M$ and if $f$ and $g$ are densities with
respect to $\mu $ such that $d\mu _0 = fd\mu $ and $d \nu _0 = gd\mu $ approach Dirac masses at
$x$ and $y$ respectively (for example by heat kernel approximations), the preceding yields
$$ P_t (Q_1 \varphi ) (x) - P_t \varphi (y) \leq  {\eu^{-2K  t} \over 2} \, d(x,y)^2, $$
that is exactly the commutation property (20). As we have seen, the latter ensures
in turn the curvature condition $CD(K  , \infty)$.

\vskip 8mm \goodbreak

{\bf 5. Alternate proofs of the commutation property}

\vskip 3mm

In this section, we briefly outline alternate proofs
of the basic commutation property (20) of Theorem 4.1. For simplicity in the notation and the exposition,
we only consider $K  =0$ below. Each proof involves at some point specific properties
and may thus be adapted to more general settings accordingly.

\vskip 2mm

{\it (i) First alternate proof.} This proof is the original argument by K. Kuwada [K1].
It requires the use of geodesics and the Hopf-Lax
formula as solution of the Hamilton-Jacobi equation. Consider,
for a smooth enough function $f : M \to \rr$,
$$ \phi (s) = P_t ( Q_s f) (x_s), \quad  s \in [0,1] ,$$
where ${(x_s)}_{s \in [0,1]}$ is a constant speed curve joining $x_0 = y$ to $x_1 = x$ in $M$.
Set for simplicity $ d = d(x,y)$. Then, by the gradient bound (6) under $CD(0 ,\infty)$,
$$ \eqalign {
    \phi' (s)  & = P_t \Big ( - {1\over 2} \, | \nabla Q_s f |^2 \Big) (x_s)
                            + \nabla P_t (Q_s f) (x_s) \cdot {\dot x}_s \cr
                    & \leq    P_t \Big ( - {1\over 2} \, | \nabla Q_s f |^2 \Big) (x_s)
                            + d \, \big | \nabla P_t (Q_s f) (x_s) \big | \cr    
                      & \leq   P_t \Big ( - {1\over 2} \, | \nabla Q_s f |^2 \Big)
                            + d \, P_t \big ( |\nabla Q_s f | \big ) \cr   
                 & \leq  {d^2 \over 2} \, . \cr } $$
Hence
$$  P_t (Q_1 f) (x) - P_t f(y) = \phi (1) - \phi (0) = \int_0^1 \phi'(s) ds \leq  {d^2 \over 2}  $$
which is the result.  

\vskip 2mm

{\it (ii) Second alternate proof.}
This second alternate proof also uses the Hopf-Lax infimum-convolution semigroup
as solution of the Hamilton-Jacobi equation,
and relies on the hypercontractivity property along the heat flow recently put forward in [B-B-G].
Namely, by the log-Harnack inequality (7), for $f$ say bounded continuous and $v >0$,
$$ P_t (Q_1f) \leq {1 \over v} \, Q_{2t} \big ( \log P_t (\eu^{v Q_1f}) \big ) . $$
Under non-negative curvature, it is shown in [B-B-G] that for every (bounded continuous)
$\psi : M \to \rr $ and $t>0$,
$$ \log P_t (e^{ Q_{2t} \psi }) \leq  P_t \psi . $$
With $ v = 1/2t$ and $ \psi = f/2t$, the conclusion immediately follows by homogeneity
of the infimum-convolutions.

\vskip 2mm

{\it (iii) Third alternate proof.} 
This proof may be obtained by linear approximations of the Hamilton-Jacobi
equation (vanishing viscosity method) along the lines of [Bo-G-L].
Following the notation therein, let for every $\varepsilon >0$,
the approximated Hopf-Lax semigroup
$$ Q_t^\varepsilon  f = - 2\varepsilon  \log P_{\varepsilon t} (\eu^{- f/2\varepsilon }) $$
solution of the equation
$$ \partial_t u = \varepsilon \, \L u  - {1\over 2} \, |\nabla u|^2 .$$
In a sense which can be made precise, $\lim_{\varepsilon \to 0} Q_t^\varepsilon  f = Q_t f$.
Dealing with
$$ \phi (s) = P_s \big (Q_1^\varepsilon (P_{t-s}f) \big), \quad  s \in [0,t],$$
shows that
$$ \phi '(s) = 2\varepsilon  P_s \bigg ( {1\over P_\varepsilon g} \,
    \bigg [ {|\nabla P_\varepsilon g|^2 \over P_\varepsilon g}
            - P_\varepsilon \bigg ( {|\nabla g|^2 \over g } \bigg ) \bigg ] \bigg) $$
where $ g = \eu^{- P_{t-s}f /2\varepsilon }$. Under the gradient bound (6),
$\phi '(s) \leq 0$ which yields that
$$  P_t(Q_1^\varepsilon f) \leq Q_1^\varepsilon (P_t f). $$
In the limit as $\varepsilon \to 0$, the announced commutation property follows.

One benefit of the third alternate proof is that it may be developed similarly
on the curvature-dimension condition $CD(0,N)$ with a finite-dimensional
parameter $N$, for example with $N=n$ for the Laplace operator
on an $n$-dimensional Riemannian manifold with non-negative Ricci curvature
(cf. [B], [Ba-G-L]). We sketch the argument.
The local logarithmic Sobolev inequalities of [B-L2] (see also [Ba-G-L])
under $CD(0,N)$ ensure after linearization that, for any $t>0$, any non-negative smooth function
$g : M \to \rr$ and any $c>0$,
$$   c \, {|\nabla P_t g|^2 \over P_t g}
            - P_t \bigg ( {|\nabla g|^2 \over g } \bigg) \leq  (c-1) P_t (\L g)  
              + {N \over 2t} \, \big (\sqrt c - 1\big )^2 P_t g . $$
Arguing as previously in the $CD(0,\infty)$ case then yields that for any $f$ and $t, s >0$,              
$$  P_t(Q_1^\varepsilon f) \leq Q_1^\varepsilon (P_s f) + N \big ( \sqrt t - \sqrt s \,\big )^2 $$
and similarly in the limit as $\varepsilon \to 0$.

Applied to the Wasserstein contraction, the latter shows that, under $CD (0,N)$
and in the notation of Corollary 4.2,
$$ W_2^2 (\mu _t , \nu _s ) \leq W_2^2 (\mu _0 , \nu_0 ) 
   + 2 N \big ( \sqrt t - \sqrt s \,\big )^2.   \leqno (22)$$
This inequality covers (21), however only when $s=t$. Note that when
$s=0$ and $\mu _t = \nu _t$, then $ W_2 (\mu _t , \mu_0) \leq \sqrt {2Nt}$
which describes a classical behavior of Brownian motion in Euclidean space.
Further Wasserstein contraction properties under curvature-dimension condition are
emphasized in [W4], [E-K-S], [K2], [B-G-G].

It should be mentioned in addition that, following the argument at the end of Section~4,
the contraction inequality (22) implies back the commutation
$$  P_t(Q_1 f) \leq Q_1 (P_s f) + N \big ( \sqrt t - \sqrt s \,\big )^2 \leqno (23) $$
for all (bounded smooth) $f : M \to \rr$ and $t,s >0$.
Now, given $a \in \rr$ to be specified, set $t = (1 + \varepsilon a) s$ for $\varepsilon >0$
small enough. By homogeneiy, the latter yields
$$  P_{(1+a\varepsilon) s}(Q_\varepsilon  f) \leq Q_\varepsilon  (P_s f) 
     + {N a^2 \varepsilon s \over  (\sqrt {1 + a\varepsilon } + 1)^2 } \, .$$
A Taylor expansion at $\varepsilon = 0 $ then shows that
$$ as P_s (\L f) - {1\over 2} \, P_s \big (|\nabla f|^2 \big ) 
    \leq - {1 \over 2} \, |\nabla P_s f|^2 + {N a^2 s \over 4} \, . $$
For $a = {2 \over N } \, \L P_s f$, it follows that
$$ |\nabla P_s f|^2 \leq P_s \big (|\nabla f|^2 \big )  - {2s \over N} \, (\L P_s f)^2 . $$
This inequality, holding (pointwise) for every (smooth) $f$ and every $s >0$, is known to be equivalent
to the curvature-dimension condition $CD (0,N)$ (cf. (12) in [B-L2], or [W4]). As such,
the Wasserstein contraction (22), as well as the dimensional commutation (23),
are also equivalent to $CD (0,N)$.

\vskip 8mm \goodbreak

{\bf 6. Links with the Evolutionary Variational Inequality}

\vskip 3mm

To conclude this work, we briefly describe some of the connections
between the preceding material and recent
contributions around the so-called Evolutionary Variational Inequality (EVI).
As mentioned in the introduction, the EVI has been recently developed by
L. Ambrosio, N. Gigli and G.~Savar\'e [A-G-S2], [A-G-S3] towards the connection
between the curvature
condition $CD(K  ,\infty)$ in the sense of the $\Gamma_2$ operator of [B-\'E] (see
[B], [Ba-G-L]), expressed here through the commutation (6), and the curvature bound in the
sense of Lott-Villani-Sturm in metric measure spaces as convexity of relative
entropy along the geodesics of optimal transportation ([L-V], [S1], [S2], [V2]).

The purpose of this short paragraph is to describe the idea at the root of the
EVI following the recent main development [A-G-S3]. In this work, the authors actually establish the
EVI in the extended class of Riemannian energy measure spaces, providing there a complete
link between the Bakry-\'Emery $\Gamma_2$ and Lott-Villani-Sturm curvatures (the implication
from Lott-Villani-Sturm to $\Gamma_2$ was achieved in [A-G-S2]).
With respect to [A-G-S3], we only concentrate here on the main principle of proof in the
simplified framework of weighted Riemannian manifolds, the main achievement
of [A-G-S3] being actually to perform the argument in a much larger class of
non-smooth spaces together with a rather involved analysis. The guideline
of this investigation is the Eulerian approach of [O-W] and [D-S]
but the non-smooth structure causes a lot of technical problems.
We nevertheless found it useful to outline the argument, avoiding all the regularity issues,
in the context of this paper to illustrate the general principle and the links
with the material of the previous sections. We of course refer to [A-G-S3]
for a complete rigorous investigation. The recent contribution [E-K-S]
addresses corresponding issues for the curvature-dimension condition $CD(K,N)$.

For simplicity thus, we deal with the weighted Riemannian framework of the preceding sections
with $d\mu  = \eu^{-V} dx$ a probability measure,
and restrict ourselves to the non-negative curvature assumption $CD(0, \infty)$ expressed
in the form of the commutation property (6) with $K  =0$. The case of arbitrary $K  \in \rr$
is easily adapted along the same lines (cf. [A-G-S3]). 

Let $f$ and $g$ be probability densities with respect to the probability measure $\mu $.
The Evolutionary Variational Inequality (EVI) indicates that under $CD(0, \infty)$, for any $t>0$,
$$ {\rm W}_2^2 (\mu _t , \nu _0 ) + 2t \int_M P_t f \log P_t f d\mu 
    \leq  {\rm W}_2^2 (\mu _0 , \nu _0 ) + 2t \int_M g \log g d \mu  \leqno (24) $$
where $d\mu _t = P_t f d\mu $, $d\nu _t = P_t g d\mu $. In the limit as $t \to 0$,
together with the semigroup property,
$$ {1 \over 2} \, {d \over dt} {\rm W}_2^2 (\mu _t , \nu _0 ) 
    \leq  \int_M g \log g d \mu  -  \int_M P_t f \log P_t f d\mu \leqno (25) $$
(the derivative being understood in an extended sense as the
limsup of the right difference quotient).

The material described in the preceding sections gets close to (24), however not quite.
Indeed, the conjunction of (12) and of the Wasserstein contraction (21) (for $K  =0$) yields 
$$ {\rm W}_2^2 (\mu _t , \nu _t ) + 2t \int_M P_t f \log P_t f d\mu 
    \leq  {3 \over 2} \, {\rm W}_2^2 (\mu _0 , \nu _0 ) + 2t \int_M g \log g d \mu   $$
which is not directly comparable to (24) but which, in any case, is useless in the limit as $t\to 0$.
To reach EVI, more on optimal transportation is actually required.

One key step in this regard is the existence of curves of probability densities $h_s$, $s \in [0,1]$,
with respect to $\mu $ interpolating between $ h_0 = g$ and $h_1 = f$,
assumed to be smooth both in space and $s$, such that
for every smooth function $\psi $ on $M$,
$$ \int_M {\dot h}_s \, \psi \,  d\mu 
    \leq   {1 \over 2}  \, {\rm W}_2^2 (\mu _0 , \nu _0 )  + {1 \over 2} \int_M |\nabla \psi |^2 h_s d\mu .
     \leqno (26) $$
Such curves are naturally provided by optimal transportation, and arise
for example in the Benamou-Brenier dynamical description of the Wasserstein
distance [B-B] (cf. [A-G-S1], [Vi1], [Vi2]). Actually,
the existence of geodesics $\mu _s = h_s \mu $
in the Wasserstein space satifying (26) just depends on the length property of the space
and is a general result of [L]. In general, it is however
not even clear that such a curve $\mu _s$ is absolutely continuous
with respect to $\mu $, so the basic issue here concerns regularity
of $\mu _s$ and $h_s$. Even in a smooth context, 
there is an correction error in (26) which may be shown to be negligible
for the further purposes so that for simplicity we ignore it here.
The existence and regularity of such curves
$h_s$, $s \in [0,1]$, satisfying (26) in a non-smooth setting is a
delicate issue carefully investigated in [A-G-S3].

To illustrate at a mild level such curves, and in
$M = \rr^n$ for the simplicity of the notation (the manifold case being
similar at the expense of further Riemannian technology, cf. [V2]),
consider the Brenier map $T : \rr^n \to \rr^n$ pushing forward
$d\mu _0 = fd\mu $ to $d\nu _0 = g d\mu $ and providing
optimal transportation in the sense of the Wasserstein distance $W_2$ as
$$  {\rm W}_2^2 (\mu _0 , \nu _0 ) =  \int_{\rrr^n} \big | x - T(x) \big |^2 f(x) d\mu (x)
   \leqno (27) $$
(cf. [A-G-S1], [V1], [V2]...).
Consider then the geodesics $T_s = s\, {\rm Id} + (1-s) T$,
$s \in [0,1]$, of optimal transportation. If $h_s $ denotes the density with respect to $\mu $
of the pushforward measure of $d \mu_0  = f d\mu $ by $T_s$
(so that $ h_0 = g$ and $h_1 = f$)
assumed to be smooth both in space and $s$, it is easily checked that
for every smooth function $\psi $ on $\rr^n$,
$$ \int_{\rrr^n} {\dot h}_s \, \psi \,  d\mu 
         = \int_{\rrr^n} \big ( x - T(x) \big ) 
             \cdot \nabla \psi  \big ( T_s (x) \big ) f(x) d\mu (x)   $$
yielding (26) by the quadratic inequality and (27). 

On the basis of (26), the EVI (24) may be analyzed by a suitable
coupling between the heat kernel and optimal transportation parametrizations.
Precisely, the expressions
$$  \int_M Q_1\varphi  P_t f d\mu  -  \int_M \varphi g d\mu 
   + t \bigg ( \int_M P_t f \log P_t f d\mu 
    -   \int_M g \log g d \mu \bigg) \leqno (28) $$
for any smooth $\varphi $ on $M$ may be represented as
$$   \int_0^1 \bigg (   {d \over ds} \int_M Q_s \varphi P_{s t} h_s d\mu 
      +   t \, {d \over ds} \int_M P_{st} h_s \log P_{st} h_s d\mu \bigg) ds. $$
Now, again under suitable smoothness assumptions not detailed here,
by the Hamilton-Jacobi equation and integration by parts,
$$  \eqalign {
     {d \over ds} \int_M Q_s \varphi P_{s t} h_s d\mu  
     & = - {1\over 2} \int_M |\nabla Q_s \varphi |^2 P_{s t} h_s d\mu 
           + \int_M {\dot h}_s P_{s t} (Q_s \varphi ) d\mu  \cr
      & \quad + t  \int_M Q_s \varphi \, {\rm L} P_{s t} h_s d\mu  \cr
     & = - {1\over 2} \int_M |\nabla Q_s \varphi |^2 P_{s t} h_s d\mu 
           + \int_M {\dot h}_s P_{s t} (Q_s \varphi ) d\mu  \cr
      & \quad - t  \int_M \nabla Q_s \varphi \cdot \nabla  P_{s t} h_s d\mu  .\cr}$$
On the other hand,
$$ \eqalign {
     {d \over ds} \int_M P_{st} h_s \log P_{st} h_s d\mu 
      & =    \int_M \big [ P_{st} {\dot h}_s + t \, {\rm L} P_{st} h_s \big ] \log P_{st} h_s d \mu \cr
      & = \int_M  P_{st} {\dot h}_s \log P_{st} h_s d\mu 
             - t \int_M  \nabla P_{st} h_s \cdot \nabla (\log  P_{st} h_s) d\mu  \cr }  $$
where we used that
$  {d \over ds}  P_{st} h_s  = P_{st} {\dot h}_s + t \, {\rm L}  P_{st} h_s  $
and $ \int_M  P_{st} {\dot h}_s d\mu  = \int_M  {\dot h}_s d\mu  = 0$.

From these expressions, it is easily checked that the sum
$$      {d \over ds} \int_M Q_s \varphi P_{s t} h_s d\mu 
      +   t \, {d \over ds} \int_M P_{st} h_s \log P_{st} h_s d\mu $$
may be rearranged as 
$$ \eqalign {
     - {1\over 2} \int_M \big |\nabla  ( Q_{s} \varphi  & 
         + t \log P_{st} h_s ) \big |^2  P_{s t} h_s d\mu 
    - {t^2 \over 2} \int_M  { |\nabla P_{st}h_s|^2 \over P_{st} h_s} \, d\mu  \cr
       & + \int_M   {\dot h}_s P_{st} \big ( Q_s \varphi  + t \log P_{st} h_s \big ) d\mu . \cr } $$
Forgetting the term
$ {t^2 \over 2} \int_M  { |\nabla P_{st}h_s|^2 \over P_{st} h_s} \, d\mu $
(which is anyway of the order $o(t)$ in the limit (25)), this
quantity is upper-bounded by
$$ - {1\over 2} \int_M P_{st} \big ( \big |\nabla  (
      Q_s \varphi   + t \log P_{st} h_s  ) \big |^2 \big )  h_s d\mu 
        + \int_M   {\dot h}_s P_{st} \big ( Q_s \varphi  + t \log P_{st} h_s \big ) d\mu $$
where we used symmetry of the semigroup. Now, by the curvature condition in the form of the
commutation (6), the latter is further upper-bounded by
$$ - {1\over 2} \int_M  
    \big | \nabla P_{st} (Q_{s} \varphi   + t \log P_{st} h_s  ) \big |^2 h_s d\mu  
        + \int_M   {\dot h}_s P_{st} \big ( Q_s \varphi  + t \log P_{st} h_s \big ) d\mu $$
With $\psi  = P_{st} (Q_s \varphi  + t \log P_{st} h_s ) $, (26) implies that this expression
is precisely bounded from above
by ${\rm W}_2^2 (\mu _0 , \nu _0 )$. Integrating in $s$ from $0$ to $1$
and taking the supremum over all $\varphi $'s then yields the
announced EVI (24). It might be worthwhile mentioning that with respect to
(6) only the weaker commutation property
$ | \nabla P_t g|^2 \leq \eu^{-K  t} P_t  ( |\nabla g|^2  ) $
is used here.   

As mentioned above, the preceding argument is inspired by the Eulerian calculus
developed by F. Otto and M. Westdickenberg [O-W] in their approach of the contraction
property (21). Namely, if the parametrization does not involve the heat flow, consider
for $\varphi  : M \to \rr$ smooth enough,
$$  \int_M Q_1\varphi  P_t f d\mu  -  \int_M \varphi P_t g d\mu 
    =   \int_0^1 \bigg (   {d \over ds} \int_M Q_s \varphi P_t h_s d\mu \bigg ) ds .$$
Since, as above,
$$  {d \over ds} \int_M Q_s \varphi P_t h_s d\mu  
      = - {1\over 2} \int_M |\nabla Q_s \varphi |^2 P_t h_s d\mu 
           + \int_M {\dot h}_s P_t (Q_s \varphi ) d\mu  , $$
by time reversibility and the gradient bound (6),
$$ \eqalign {
    {d \over ds} \int_M Q_s \varphi P_t h_s d\mu  
      & = - {1\over 2} \int_M P_t \big (|\nabla Q_s \varphi |^2)  h_s d\mu 
           + \int_M {\dot h}_s P_t (Q_s \varphi ) d\mu  \cr
      & \leq  - {1\over 2} \int_M \big | \nabla P_t ( Q_s \varphi) \big |^2  h_s d\mu 
           + \int_M {\dot h}_s P_t (Q_s \varphi ) d\mu .  \cr }$$
Using (26) then yields
$$  \int_M Q_1\varphi  P_t f d\mu  -  \int_M \varphi P_t g d\mu 
     \leq  {1 \over 2} \, W_2^2 ( \mu _0, \nu _0) ,$$
that is, after taking the supremum in $\varphi $, the contraction property (21)
of Corollary 4.2.

\vskip 8mm

{\it Acknowledgement. We are thankful to L. Ambrosio, N. Gigli and G. Savar\'e for helpful
discussions on the EVI and for pointing out relevant references, and to A. Guillin for
sharing with us his simple proof
from the Wasserstein contraction to the curvature condition at the end of Section 4.
We are also most grateful to the referee for numerous comments and corrections that helped
in improving the exposition.}

\vskip 12mm

\goodbreak

\baselineskip=10pt

\font\pc=cmcsc10 \rm 

\centerline {\pc References}

\vskip 5mm

\font\pc=cmcsc8 \rm 

\font\eightrm =cmr8 

{\eightrm 

\ref A-G-S1|{\pc L. Ambrosio, N. Gigli, G. Savar\'e.}
Gradient flows in metric spaces and in the
space of probability measures. Lectures in Mathematics ETH Z\"urich. Birkh\"auser (2008)|

\ref A-G-S2|{\pc L. Ambrosio, N. Gigli, G. Savar\'e.} Metric measures spaces with Riemannian
Ricci curvature bounded from below (2011)|

\ref A-G-S3|{\pc L. Ambrosio, N. Gigli, G. Savar\'e.} Bakry-\'Emery curvature-dimension
condition and Riemannian Ricci curvature bounds (2012)|

\ref B|{\pc D. Bakry.} L'hypercontractivit\'e et son utilisation en th\'eorie des semigroupes.
Ecole d'Et\'e de Probabilit\'es de Saint-Flour. Lecture Notes in Math. 1581, 1--114 (1994).
Springer|

\ref B-B-G|{\pc D. Bakry, F. Bolley, I. Gentil.} Dimension dependent 
hypercontractivity for Gaussian kernels. Probab. Theory Related Fields 154, 845--874 (2012)|

\ref B-\'E|{\pc D. Bakry, M. \'Emery.} Diffusions hypercontractives.
S\'eminaire de Probabilit\'es {XIX}, 1983/84, Lecture Notes in Math. 1123, 177--206. Springer (1985)|

\ref Ba-G-L|{\pc D. Bakry, I. Gentil, M. Ledoux.} Analysis and geometry of Markov
diffusion operators. Grundlehren der Mathematischen Wissenschaften 348. Springer (2013)|

\ref B-L1|{\pc D. Bakry, M. Ledoux.} L\'evy-Gromov's isoperimetric inequality for an
infinite-dimensional diffusion generator. Invent. Math. 123, 259--281 (1996)|

\ref B-L2|{\pc D. Bakry, M. Ledoux.} A logarithmic Sobolev form of the Li-Yau parabolic
inequality. Rev. Mat. Iberoam. 22, 683--702 (2006)|

\ref B-B|{\pc J.-D. Benamou, Y. Brenier.} A computional fluid mechanics solution to the
Monge-Kantorovich mass transfer problem. Numer. Math. 84, 375--393 (2000)|

\ref Bo-G-L|{\pc S. Bobkov, I. Gentil, M. Ledoux.} Hypercontractivity of Hamilton-Jacobi equations.
J. Math. Pures Appl. 80, 669--696 (2001)|

\ref B-H|{\pc S. Bobkov, C. Houdr\'e.} Some connections between isoperimetric and Sobolev-type
inequalities. Mem. Amer. Math. Soc. 129 (1997)|

\ref B-G-G|{\pc F. Bolley, I. Gentil, A. Guillin.} Dimensional contraction via Markov
transportation distance (2013)|

\ref C-MC-V|{\pc J. Carrillo, R. McCann, C. Villani.}
Kinetic equilibration rates for granular media and related
equations: entropy dissipation and mass transportation estimates.
Rev. Mat. Iberoam. 19, 971--1018 (2003)|

\ref CE|{\pc D. Cordero-Erausquin.} Some applications of mass transport to Gaussian-type
inequalities. Arch. Ration. Mech. Anal. 161, 257--269 (2002)|

\ref D-S|{\pc S. Daneri, G. Savar\'e.} Eulerian calculus for the displacement convexity
in the Wasserstein distance. SIAM J. Math. Anal. 40, 1104-1122 (2008)|

\ref D|{\pc E. B. Davies.} Heat kernels and spectral theory.
Cambridge Tracts in Mathematics 92. Cambridge (1989)|

\ref Er|{\pc M. Erbar.} The heat equation on manifolds as a gradient flow
in the Wasserstein space. Ann. Inst. H. Poincar\'e -- Probabilit\'es et Statistiques 46, 1--23 (2010)|

\ref E-K-S|{\pc M. Erbar, K. Kuwada, K.-T. Sturm.} On the equivalence of the entropic 
curvature-dimension condition and BochnerÕs
inequality on metric measure spaces (2013)|

\ref Ev|{\pc L. Evans.} Partial differential equations. Graduate Studies in Mathematics 19.
American Mathematical Society (1998)|

\ref G-K-O|{\pc N. Gigli, K. Kuwada, S.-i. Ohta.} Heat flow on Alexandrov spaces.
Comm. Pure Appl. Math., to appear (2012)|

\ref G-L|{\pc N. Gigli, M. Ledoux.} From log Sobolev to Talagrand: a quick proof. 
Discrete and Continuous Dynamical Systems 33, 1927--1935 (2013)|

\ref H|{\pc M. Hino.} On short time asymptotic behavior of some symmetric
diffusions on general state spaces. Potential Anal. 16, 249--264 (2002)|

\ref J-K-L|{\pc R. Jordan, D. Kinderlehrer, F. Otto.} The variational formulation of
the Fokker-Planck equation. SIAM J. Math. Anal. 29, 1--17 (1998)|

\vfill
\supereject

\ref K1|{\pc K. Kuwada.} Duality on gradient estimates and Wassertein controls.
J. Funct. Anal. 258, 3758--3774 (2010)|

\ref K2|{\pc K. Kuwada.} Space-time Wasserstein controls and 
Bakry-Ledoux type gradient estimates (2013)|

\ref L|{\pc S. Lisini.} Characterization of absolutely continuous curves in
Wasserstein spaces. Calc. Var. Partial Differential Equations 28, 85--120 (2007)|

\ref L-Y|{\pc P. Li, S.-T. Yau.} On the parabolic kernel of the Schr\"odinger operator.
Acta Math. 156, 153--201 (1986)|

\ref L-V|{\pc J. Lott, C. Villani.} Ricci curvature for metric-measure spaces via optimal transport.
Ann. Math. 169, 903--991 (2009)|

\ref O|{\pc F. Otto.} The geometry of dissipative evolution equations: the porous
medium equation. Comm. Partial Differential Equations 26, 101--174 (2001)|

\ref O-V|{\pc F. Otto, C. Villani}. Generalization of an inequality by Talagrand,
and links with the logarithmic Sobolev inequality. J. Funct. Anal. 173, 361--400 (2000)|

\ref O-W|{\pc F. Otto, M. Westdickenberg.} Eulerian calculus for the contraction in the
Wasserstein distance. SIAM J. Math. Anal. 37, 1227--1255 (2005)|

\ref vR-S|{\pc M.-K. von Renesse, K.-T. Sturm}.
Transport inequalities, gradient estimates, entropy and Ricci Curvature.
Comm. Pure Appl. Math. 68, 923--940 (2005)|

\ref S1|{\pc K.-T. Sturm}. On the geometry of metric measure spaces I. Acta Math. 196, 65--131 (2006)|

\ref S2|{\pc K.-T. Sturm}. On the geometry of metric measure spaces II. Acta Math. 196, 133-177 (2006)|

\ref V1|{\pc C. Villani.} Topics in optimal transportation. Graduate
Studies in Mathematics 58. American Mathematical Society (2003)|

\ref V2|{\pc C. Villani.} Optimal transport. Old and new.
Grundlehren der Mathematischen Wissenschaften 338. Springer (2009)|

\ref W1|{\pc F.-Y. Wang.} Logarithmic Sobolev inequalities on noncompact Riemannian manifolds.
Probab. Theory Related Fields 109, 417--424 (1997)|

\ref W2|{\pc  F.-Y. Wang.} Functional inequalities, Markov properties and spectral theory.
Science Press (2005)|

\ref W3|{\pc  F.-Y. Wang.} Harnack inequalities on manifolds with boundary and
applications. J. Math. Pures Appl. 94, 304--321 (2010)|

\ref W4|{\pc  F.-Y. Wang.} Equivalent semigroup properties for the curvature-dimension condition.
Bull. Sci. Math. 135, 803--815 (2011)|

}

\vskip 10mm

\parskip 1,2mm

\font\eightrm =cmr8  {\eightrm

{\pc D. B., M. L.: Institut de Math\'ematiques de Toulouse, Universit\'e de Toulouse,
31062 Toulouse, France, and Institut Universitaire de France,   
\font\eighttt =cmtt8  {\eighttt  bakry, ledoux@math.univ-toulouse.fr}} 

{\pc I. G.:  Laboratoire Camille Jordan, Universit\'e de Lyon, 69622 Lyon, France,   
\font\eighttt =cmtt8  {\eighttt  gentil@math.univ-lyon1.fr}}

}

\bye